\documentclass[11pt]{article}
\usepackage{authblk} 

\usepackage[latin1]{inputenc}
\usepackage[T1]{fontenc}
\usepackage{amsmath}
\usepackage{amsfonts}
\usepackage{amssymb}
\usepackage{enumerate}
\usepackage[normalem]{ulem} 

\usepackage{amsmath,xcolor,ulem}
\usepackage{amsfonts}
\usepackage{amssymb,todonotes}
\usepackage{amsthm}
\usepackage{graphicx}
\usepackage{epsfig}
\usepackage{subcaption}
\usepackage{marginnote}

\usepackage[colorlinks=true, allcolors=blue]{hyperref}

\usepackage[top=2.8cm,bottom=2.8cm,left=2.5cm,right=2.5cm]{geometry}
\usepackage{float}
\usepackage[font=footnotesize,width=\linewidth]{caption}

\def\H{\mathcal H}
\def\U{\mathcal U}
\def\d{\mathrm d}
\def\R{\mathbb R}

\newcommand{\iu}{\mathrm{i}} 
\newcommand*\xbar[1]{%
  \,\hbox{%
    \vbox{%
      \hrule height 0.1pt 
      \kern0.4ex
      \hbox{%
        \kern-0.1em
        \ensuremath{#1}%
        \kern-0.1em
      }%
    }%
  }\,%
} 

\title{Modelling car-following dynamics with stochastic input-state-output port-Hamiltonian systems}
\author{Julia Ackermann$^1$, Matthias Ehrhardt$^1$, Thomas Kruse$^1$ and Antoine Tordeux$^2$}
\affil{$^1$Applied and Computational Mathematics, University of Wuppertal, Germany}
\affil{$^2$Traffic Safety and Reliability, University of Wuppertal, Germany}

\begin{document}
\maketitle

\begin{tikzpicture}[remember picture,overlay]
	\node[anchor=north east,inner sep=20pt] at (current page.north east)
	{\includegraphics[scale=0.2]{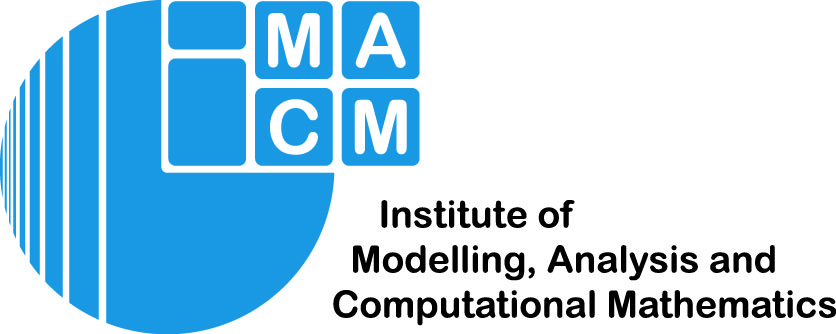}};
\end{tikzpicture}

\begin{abstract}
In this contribution, we introduce a general class of car-following models with an input-state-output port-Hamiltonian structure. 
 We derive stability conditions and long-term behavior of the finite system with periodic boundaries and quadratic interaction potential by spectral analysis and using asymptotic properties of multivariate Ornstein-Uhlenbeck processes. 
 The uncontrolled dynamics exhibit instability and random collective behavior under stochastic perturbations. 
 By implementing an open-loop speed control, the system stabilizes and weakly converges to Gaussian limit distributions. 
 The convergence is unconditional for constant speed control.
 However, a stability condition arises for the closed-loop system where the speed control acts as a dynamic feedback depending on the distance ahead. 
 The results are illustrated by numerical simulations. 
 Interestingly, only the closed-loop system is able to reproduce, at least transiently, realistic stop-and-go behavior that can be resolved using the Hamiltonian component of the model.
\end{abstract}

\begin{minipage}{0.9\linewidth}
\medskip

\end{minipage}

\section{Introduction}
Port-Hamiltonian systems (PHS) have become a fundamental framework for modeling and analyzing complex physical systems 
due to their ability to incorporate energy-based structures and dissipation phenomena. 
This mathematical paradigm has been widely applied in domains such as robotics, electrical circuits, and thermodynamics, where energy conservation and exchange are central to system dynamics  \cite{rashad2020twenty}.
Inspired by recent results of Knorn et al.\ \cite{knorn2014passivity} and Matei et al.\ \cite{matei2019}, 
we have formulated a general class of stochastic car-following models with an input-state-output port-Hamiltonian structure \cite{ehrhardt2024collective,ackermann2024stabilisation,rudiger2024stability}. 
This port-Hamiltonian formulation introduces a new modeling paradigm for the analysis of car-following dynamics. 

In this work, we present a novel class of stochastic car-following models formulated within the input-state-output port-Hamiltonian framework.
Our approach builds on earlier studies, particularly the 
car-following dynamics explored by Bando et al.\ \cite{Bando1995} and Jiang et al.\ \cite{jiang2001full} based on \emph{optimal velocity functions}. 
However, in contrast to classical approaches where the interaction is totally asymmetric, the port-Hamiltonian car-following models include both vehicles in front and behind in the interaction. 
By leveraging the PHS formulation, we systematically incorporate noise and external control inputs into the car-following dynamics,
thereby enabling a comprehensive analysis of both uncontrolled and controlled traffic systems. 
This extension provides new insights into stability properties and emergent behaviors such as stop-and-go waves, 
which are currently observed in real-world traffic and in experiments \cite{Orosz2009,Sugiyama2008}.

The proposed models do not claim to provide a more realistic description of driver behavior. 
However, they extend previous research by introducing a rigorous mathematical formulation of car-following dynamics as multivariate stochastic differential equations with periodic boundary conditions. 
Notably, the Hamiltonian framework enables us to analyze the asymptotic behavior of these systems, 
leveraging spectral theory and properties of multivariate Ornstein-Uhlenbeck (OU) processes. 
The stochastic nature of the model captures realistic fluctuations in driver behavior, 
such as delays and variability in acceleration, which are critical for reproducing observed traffic phenomena.

In uncontrolled systems, the absence of active stabilization mechanisms leads to instability and collective random behavior characterized by the divergence of the ensemble's mean velocity. 
By contrast, open-loop control stabilizes the system, ensuring weak convergence to Gaussian limit distributions under constant speed control. 
Moreover, closed-loop systems, where the speed control dynamically depends on the distance to the preceding vehicle, 
exhibit richer behaviors, including transient stop-and-go patterns. 
These observations align with empirical findings and highlight the robustness of the PHS-based approach to modeling complex traffic dynamics.

Our study is closely related to recent advancements in adaptive cruise control systems and their mathematical underpinnings \cite{ISO15822, kesting2008adaptive}. 
We show that the Hamiltonian structure not only simplifies the analysis of stability conditions but also enables the systematic design of control laws to mitigate traffic instabilities.
In particular, we derive exact stability criteria for closed-loop systems, demonstrating that increasing the Hamiltonian component or control parameters improves system stability.
These results are complemented by numerical simulations, which illustrate key behaviors and validate the theoretical predictions.

The remainder of this paper is organized as follows.
In Section~\ref{sec:2}, we introduce the stochastic port-Hamiltonian car-following model and outline its mathematical formulation. 
Section~\ref{sec:3} analyzes the asymptotic behavior of uncontrolled, open-loop, and closed-loop systems, with a focus on stability and convergence properties. 
Numerical simulations in Section~\ref{sec:num} validate the theoretical findings and provide insights into transient and long-term behaviors.

\section{Stochastic port-Hamiltonian car-following model}\label{sec:2}

\subsection{Notations}
Consider $N$ vehicles on a segment of length $L$ with periodic boundaries. 
The positions of the $N$ vehicles at time $t\ge0$ are denoted by $q_1(t),\dots,q_N(t)\in\R$ while their speeds are the variables $p_1(t),\dots,p_N(t)\in\R$. 
Assume that the vehicles are initially ordered by their indices, i.e.,
$0 \le q_1(0) < q_2(0) < \ldots < q_N(0)\le L$.
The distances and the speed differences to the preceding vehicle are given by
\begin{equation}
    \Delta q_n(t)=\left\{\begin{array}{lll}
         q_{n+1}(t)-q_n(t),&\qquad\quad& n<N,\\
         L+q_1(t)-q_N(t),&&n=N,
    \end{array}\right.
\end{equation}
and
\begin{equation}
    \Delta p_n(t)=\left\{\begin{array}{lll}
         p_{n+1}(t)-p_n(t),&\qquad\quad& n<N,\\
         p_1(t)-p_N(t),&&n=N,
    \end{array}\right.
\end{equation}
respectively. 

\subsection{Car-following model}
A general symmetric car-following model in one dimension is given by
\begin{equation}
\left\{\begin{aligned}\label{Eq:ModMicro}
         ~\d q_n(t) &= p_n(t)\,\d t,\\[.5mm]
         ~\d p_n(t) &= \Bigl[\gamma \bigl[u_n(t)-p_n(t)\bigr]
           + \beta\bigl[\Delta p_n(t)-\Delta p_{n-1}(t)\bigr]\\[-.5mm]
         &\qquad\quad+\U'\bigl(\Delta q_n(t)\bigr)-\U'\bigl(\Delta q_{n-1}(t)\bigr)\Bigr]\,\d t +  \sigma \,\d W_n(t),
\end{aligned}\right.
\end{equation}
where 
\begin{itemize}
    \item $u_n$ is the input component that acts as a speed control with relaxation rate $\gamma\ge0$.
    It is an external control parameter for open-loop systems, or can be state-dependent and acts as a feedback control for closed-loop systems.
    \item $\beta\ge0$ is the speed alignment rate with the nearest neighbors and corresponds to the dissipative part of the PHS. 
    \item $\U'$ is the derivative of a distance-based (isotropic) convex interaction potential $\U\in C^1(\R,[0,\infty))$, which allows collisions to be avoided.
    It provides the potential energy in the Hamiltonian and, together with the kinematic relations, the skew-symmetric structure of the PHS. 
    \item $W=\bigl(W_n\bigr)_{n=1}^N\colon [0,\infty)\times \Omega \to \R^N$ is an $N$-dimensional standard Brownian motion on a probability space $(\Omega, \mathcal{F}, \mathbb{P})$ that models noise in the dynamics with volatility $\sigma \in \R\backslash\{0\}$. 
\end{itemize}

\subsection{Port-Hamiltonian formulation}
With $z=(\Delta q,p)^\top$, the port-Hamiltonian formulation is given by the stochastic input-state-output system 
\begin{equation}\label{Eq:PHS}
    \d z(t) = \bigl(J-R\bigr)\nabla \H\bigl(z(t)\bigr)\,\d t
               + U\bigl(z(t)\bigr)\,\d t
                   + \Sigma\,\d W(t),
\end{equation}
where the Hamiltonian operator is given by
\begin{equation}
    \H\binom{\Delta q}{p}=\frac{1}{2}\sum_{n=1}^N p_n^\top
    p_n+\sum_{n=1}^N\U(\Delta q_n),
\end{equation}
and where
\begin{equation*}
A=\begin{bmatrix}-1&1&\\[-1mm]
  &\ddots&\ddots\\[0mm]
  1&&-1~\end{bmatrix}\in \R^{N\times N}, \quad
  J=\begin{bmatrix}
        0&A\\-A^\top&0
    \end{bmatrix}\in \R^{2N\times 2N}
    \;\text{skew symmetric},
    \end{equation*}
 \begin{equation*}   
 U(z)=\bigl[0,\gamma u(z)\bigr]^\top\in\R^{2N},\quad 
  R=\begin{bmatrix}
        0&0\\0&\beta A^\top A+\gamma I_N
    \end{bmatrix}\in \R^{2N\times 2N}
    \;\text{positive semidefinite},
    \end{equation*}
  $u\colon\R^{2N}\to \R^{N}$ is the input control,
  $\Sigma=\bigl[0, \sigma I_N\bigr]^\top\in \R^{2N\times N}$ is the noise volatility.

The study focuses on the control terms $\gamma$ and $u$, which seem to play a dominant role in the dynamics. 
We successively consider the uncontrolled case where $\gamma=0$ before considering the cases $\gamma>0$ where, in an open-loop framework, 
$u_n$ is constant and equal for all the agents and where, in a closed-loop framework,
$u_n$ is a feedback controller that depends linearly on the gap between the vehicles.

\section{Long-term behavior}\label{sec:3}

\subsection{Uncontrolled system}
Without control, i.e., when $\gamma=0$, the system is Hamiltonian dissipative
\begin{equation}
    \d z(t) = \bigl(J-R\bigr)\nabla \H\bigl(z(t)\bigr)\,\d t + \Sigma\,\d W(t).
\end{equation}
The noise plays the role of the inputs in the port-Hamiltonian formulation.
In this case, the stochastic system does not converge to a limiting distribution \cite{ehrhardt2024collective}. 
In fact, thanks to the {telescopic form} of the model and the periodic boundaries, 
we can quickly observe that the ensemble's mean velocity
\begin{equation}\label{eq:average_velocity}
    \xbar{p}(t)
    =\frac{1}{N}\sum_{n=1}^N p_n(t), \quad t\in [0,\infty),
\end{equation}
is a diverging Brownian motion 
for any potential $\U\in C^1(\R,[0,\infty))$. 
For $\alpha\in (0,\infty)$, we consider in the sequel the quadratic distance-based potential $\U(x)=(\alpha x)^2/2$. 
The Hamiltonian and its gradient then are
\begin{equation*}
\H(z)=\frac{1}{2} \|p\|^2+\frac{\alpha^2}{2}\|\Delta q\|^2 \quad\text{and}
\quad\nabla \H(z)=\bigl[\alpha^2 \Delta q,p\bigr]^\top,
\end{equation*}
while the system is the multidimensional Ornstein-Uhlenbeck stochastic process
\begin{equation}\label{eq:OU}
    \begin{aligned}
        \d z(t) = B z(t)\,\d t + \Sigma\,\d W(t),
    \end{aligned}
\end{equation}
where 
\begin{equation}
   B = \begin{bmatrix} 0&A\\ -\alpha^2A^\top &-\beta A^\top A \end{bmatrix}
    \in \R^{2N\times 2N}\quad\text{and}\quad
     \Sigma = \begin{bmatrix}0\\\sigma I_N\end{bmatrix}\in\R^{2N\times N}.
\end{equation}
Each block $N\times N$ of $B$ is circulant and we have for $\lambda\in\mathbb C$ the determinant 
\begin{equation}
        |B-\lambda I_{2N}|=|C|\quad\text{with}\quad C=\lambda^2 I_{N}
            + \lambda \beta A^\top A + \alpha^2 A^\top A.
\end{equation}
The matrix $C\in \R^{N\times N}$ is also circulant. Since the determinant of a circulant matrix 
\begin{equation*}
C=\begin{bmatrix}
    c_0 & c_{N-1} & \ldots & c_1\\
    c_1 & c_0 & \ldots & c_2\\
    \vdots & \vdots & \ddots &\vdots \\
    c_{N-1} & c_{N-2} & \ldots & c_0
    \end{bmatrix}
\end{equation*}
is the product
\begin{equation*}
    |C| = \prod_{j=0}^{N-1} \bigl[c_0 + c_{N-1} \omega^j + c_{N-2} \omega^{2j} + \dots + c_1\omega^{(N-1)j}\bigr],
\end{equation*}
where $\omega=\exp\big({2\pi\iu/N}\big)$,
the characteristic equation of the matrix $B$ is
\begin{equation}
    \prod_{j=0}^{N-1} \bigl[\lambda^2 + \lambda \beta\mu_j + \alpha^2\mu_j\bigr]=0,\quad\text{with}\quad
    \mu_j=2-2\cos\biggl(\frac{2\pi j}{N}\biggr)\in [0,4].
    \end{equation}
From this we can directly compute the eigenvalues of $B$.
The real parts of the $2(N-1)$ eigenvalues for $j=1,\ldots,N-1$ are non-positive. 
However, the two zero eigenvalues for $j=0$ make the process $(z(t))_{t\ge 0}$ diverging \cite{ehrhardt2024collective}.

The deviation from the ensemble's mean velocity 
\begin{equation}
   p_n(t)-\xbar{p}(t)
   =\Bigl(1-\frac{1}{N}\Bigr) p_n(t) - \frac{1}{N}\sum_{k\neq n}p_k(t),\qquad t\ge 0,\quad  n =1,\ldots,N,
\end{equation}
gives rise to 
the $\R^{2N}$-valued process $(x(t))_{t\ge 0}$ given by  $x(t)
=\bigl(\Delta q(t), M p(t)\bigr)^\top$
 where $M \in \R^{N\times N}$ is the deviation from the ensemble's mean velocity metric. It 
follows the dynamics
\begin{equation}
    \d x(t)=Bx(t)\,\d t
     +\begin{bmatrix}0\\ \sigma M \end{bmatrix}\,\d W(t).
\end{equation}
Here, the noise averaging $\sigma M$ causes the deviation process $(x(t))_{t\ge 0}$ to converge \cite{ehrhardt2024collective}. 
In summary, the mean speed of the uncontrolled vehicles diverges as a Brownian motion. Interestingly, however, the variance of the vehicle speeds converges.

\subsection{Open-loop system}
With a constant input speed control, i.e., $\gamma>0$ and $u_n\equiv x\in\R$ for all $n\in\{1,\ldots,N\}$, the system is a stochastic input-state-output port-Hamiltonian system. 
In contrast to the previous model without control, the ensemble's mean velocity \eqref{eq:average_velocity} converges.
In fact, thanks to the model's telescopic form and periodic boundaries, the mean velocity follows an Ornstein-Uhlenbeck process for any potential function $\U\in C^1(\R,[0,\infty))$:
\begin{equation}\begin{aligned}
   \d\xbar{p}(t)
   = \gamma \bigl(x-\xbar{p}(t)\bigr)\,\d t
     +\frac{\sigma}N\sum_{n=1}^N \d W_n(t), \quad t\in [0,\infty).
\end{aligned}
\end{equation}
In addition, the ensemble's mean velocity \eqref{eq:average_velocity} relaxes to $x$ for the deterministic system with $\sigma=0$.

We consider the dynamics of the difference to the controlled velocity $x$ to deal with a linear system
\begin{equation}
    \tilde{p}_n(t)=p_n(t)-x,\quad t\in[0,\infty),\; n\in\{1,\ldots,N\}.
\end{equation} 
The process $(\Delta q,\tilde{p})^\top$ with 
the quadratic potential $\U(x)=(\alpha x)^2/2$ is the multidimensional Ornstein-Uhlenbeck process \eqref{eq:OU} where
\begin{equation}
   B=\begin{bmatrix}0&A\\-\alpha^2A^\top &-\beta A^\top A -\gamma I_N\end{bmatrix}
    \in \R^{2N\times 2N}\quad\text{and}\quad
     \Sigma=\begin{bmatrix}0\\\sigma I_N\end{bmatrix}\in\R^{2N\times N}.
\end{equation}
Each block $N\times N$ of $B$ is circulant, and we obtain the characteristic equation of the open-loop system
\begin{equation}
       \prod_{j=0}^{N-1} \bigl[ \lambda^2 + \lambda \bigl(\beta\mu_j + \gamma\bigr) + \alpha^2\mu_j\bigr]=0,\quad \mu_j=2-2\cos\left(2\pi j/N \right)\in [0,4].
    \end{equation}
From this we can directly compute the eigenvalues of $B$. 
The system is truncated due to the periodic boundaries (the sum of the distances is systematically $L$ at all times) and one eigenvalue for $j=0$ is zero. 
However, all the remaining eigenvalues $\lambda$ have 
negative 
real parts. 
This makes the open-loop system unconditionally stable and convergent \cite{ackermann2024stabilisation}.

\subsection{Closed-loop system}
More realistic models, related to optimal speed models in traffic engineering \cite{Bando1995, jiang2001full},  
have a feedback control 
$u_n \equiv f(\Delta q_n)$ with $f \in C^1(\R,\R)$ 
where the control speed depends on the distance ahead. 
They can describe complex instability phenomena such as stop-and-go waves \cite{Orosz2009, rudiger2024stability}.
In addition, such a modeling component allows to include a constant time gap parameter $T>0$ for linear functions $f$, which is the following strategy observed in real traffic, recommended by international standards \cite{ISO15822}, and used to design adaptive cruise control systems \cite{kesting2008adaptive}.
The state-dependent input control plays the role of a feedback in the dynamics. 
However, in contrast to the constant speed control \cite{ackermann2024stabilisation}, 
the system is no longer unconditionally stable. 
The stochastic car-following model is 
\begin{equation}\label{Eq:ModMicro2}
\left\{\begin{aligned}
         ~\d q_n(t)&=p_n(t)\,\d t,\\[.5mm]
         ~\d p_n(t)&=\Bigl[{\gamma\bigl[f\bigl(\Delta q_n(t)\bigr)-p_n(t)\bigr]}+
         {\beta\bigl[\Delta p_n(t)-\Delta p_{n-1}(t)\bigr]}\\[-1mm]
         &\qquad\quad+{\U'\bigl(\Delta q_n(t)\bigr)-\U'\bigl(\Delta q_{n-1}(t)\bigr)}\Bigr]\d t 
         +  {\sigma \,\d W_n(t)}.
\end{aligned}\right.
\end{equation}
Consider the quadratic interaction potential $\U(x)=(\alpha x)^2/2$, the linear feedback $f(x)=(x-\ell)/T$, with $T>0$ the time gap and $\ell\ge0$ the vehicle size, and the perturbed system with 
\begin{equation*}
    \tilde p_n(t)=p_n(t) + \ell/T, \quad t\in[0,\infty), ~~n\in\{1,\ldots,N\}.
\end{equation*}
Then the closed-loop system is the multidimensional Ornstein-Uhlenbeck process \eqref{eq:OU} 
where
\begin{equation}
   B=\begin{bmatrix}0&A\\\frac\gamma T I_N-\alpha^2A^\top &-\beta A^\top A -\gamma I_N\end{bmatrix}
    \in \R^{2N\times 2N}\quad\text{and}\quad
     \Sigma=\begin{bmatrix}0\\\sigma I_N\end{bmatrix}\in\R^{2N\times N}.
\end{equation}\vspace{0mm}
Here again, each block $N\times N$ of $B$ is circulant and we obtain the characteristic equation
\begin{equation}
       \prod_{j=0}^{N-1} \bigl[\lambda^2 + \lambda (\beta\mu_j + \gamma) + \alpha^2\mu_j+\frac\gamma T(1-\omega^j)\bigr]=0,\quad~~ 
			\left|\begin{array}{l}
			\mu_j=2-2\cos\left(2\pi j/N\right),\\[1mm]
			\omega=\exp(2\pi\iu/N) .
			\end{array}\right.
    \end{equation}
From this we can directly compute the eigenvalues of $B$: $\lambda_{0,0}=0$ and $\lambda_{0,1}=-\gamma$ for $j=0$ and 
\begin{equation}
\lambda_{j,k}=\frac{1}{2}\bigg[-(\beta\mu_j+\gamma)+(-1)^k\sqrt{(\beta\mu_j+\gamma)^2-4\Bigr(\alpha^2 \mu_j+\frac\gamma T(1-\omega^j)\Bigr)}\,\bigg]
\end{equation}
for all $j\in \{1,\ldots, N-1\}$ and $k\in \{0,1\}$.
    
For simplicity, we introduce 
$c_j = \cos(2\pi j/N)$, 
$s_j =\sin(2\pi j/N)$ and 
\begin{equation}
	\left|~\begin{aligned}
		\kappa_j &= 2\beta(1-c_j) + \gamma,\\      
		\eta_j &= 0,                              
	\end{aligned}\right. \qquad 
	\left|~\begin{aligned}
		\nu_j &= (1-c_j) \bigl(\textstyle\frac\gamma T+2\alpha^2\bigr),\\
		\rho_j &=  -\textstyle\frac\gamma T s_j.
	\end{aligned}\right.
\end{equation}
Recall that the roots of the 2nd order polynomial $\lambda_j^{2}+(\kappa_j+\iu \eta_j) \lambda_j +\nu_j+\mathrm{i} \rho_j=0$ satisfy $\operatorname{Re}(\lambda_j) <0$ if and only if \cite[Th.~3.2]{Frank1946}
\begin{equation}
     \kappa_j>0\quad\text{and}\quad 
     \det\begin{bmatrix}
		\kappa_j & 0 & -\rho_j \\
		1 & \nu_j & -\eta_j \\
		0 & \rho_j & \kappa_j
	\end{bmatrix}
    =\kappa_j\left(\nu_j \kappa_j+\rho_j \eta_j\right)-\rho_j^{2}>0.
\end{equation}
We obtain the exact stability conditions
\begin{equation}
	\begin{aligned}
        \gamma>0,\quad\bigl(2\beta(1-c_j)+\gamma\bigr)^2\,
        \Bigl(\frac{\gamma}{T}+2 \alpha^2\Bigr) 
        -\Bigl(\frac{\gamma}{T}\Bigr)^2 \bigl(1+c_j\bigr) >0,\quad j\in \{1,\ldots, N-1\},
	\end{aligned}
\end{equation}
which simplify to the sufficient conditions 
\begin{equation}\label{eq:CS1}
    \gamma>0,\quad~~ \gamma T+2(\alpha T)^2 > 2.
\end{equation}
Increasing the Hamiltonian component $\alpha$, the time gap $T$, or the control rate $\gamma$ improves the stability.
In addition, the speed alignment parameter has only a small effect.

\section{Numerical simulation}\label{sec:num}

We simulate $N=20$ vehicles on a segment of length $L=141$ with periodic boundaries, starting from uniform initial conditions\footnote{An online simulation module of a bigger system (to make the motion online faster) is available at the address\\ \url{https://www.vzu.uni-wuppertal.de/fileadmin/site/vzu/Port-Hamiltonian_single-file_models.html?speed=0.8}.}. 
The simulations are performed using an implicit Euler scheme for the positions of the vehicles and an explicit Euler-Maruyama scheme for the vehicle speeds. 
The time step is $\delta t=0.001$ for both numerical schemes.
We repeat three independent simulations for each of the uncontrolled system where $\gamma=0$, the open-loop system where $\gamma=0.1$ and $u_n=x=(L/N-\ell)/T=2.05$ for all the vehicles $n\in\{1,\ldots,N\}$, and the closed-loop system where $\gamma=1$ and $u_n=:(\Delta x_n-\ell)/T$ with $\ell=5$ and $T=1$. 
Note that the controlled speed $x=(L/N-\ell)/T=2.05$ in the open-loop system is set to match the speed of uniform solutions for the closed-loop system. 
We systematically show the trajectories of the vehicles and the vehicle mean speed \eqref{eq:average_velocity}, the speed of a single vehicle, and the empirical speed variance 
\begin{equation}
    V(t)=\frac1{N-1}\sum_{n=1}^N\bigl(p_n(t)-\xbar p(t)\bigr)^2
\end{equation}
over the first 250 time units in an adjacent plot.

\begin{figure}[H]
    \centering
    \includegraphics[width=0.87\linewidth]{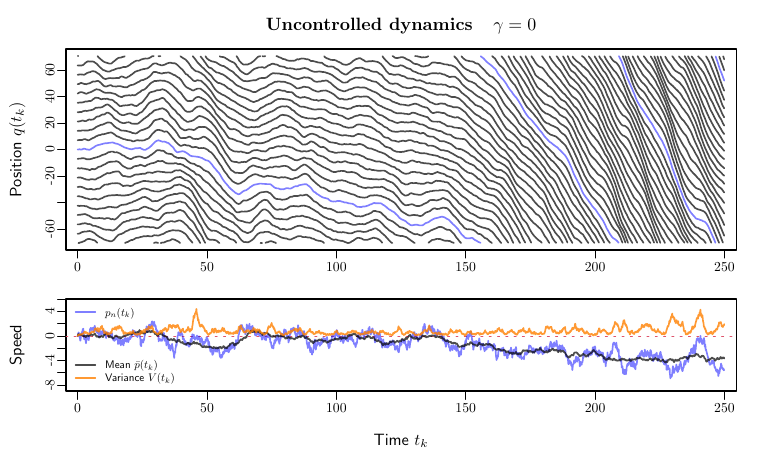}
    \caption{Vehicle trajectories (upper panel) and vehicle mean speed, speed variance and speed of a single vehicle (lower panel) for the uncontrolled system with $\alpha=1$, $\beta=1$, $\gamma=0$; $\sigma=1$. The mean speed diverges as a Brownian motion while the speed variance converges. This makes the vehicles collectively move in random directions.}
    \label{fig1}
\end{figure}

The simulation of the uncontrolled system is shown in Figure~\ref{fig1}.
The mean speed of the vehicles is a Brownian motion reaching some negative values in this example. 
Since the Brownian motion diverges but the speed variance converges, after a while we can expect all agents to move together in a given direction,
even if the interaction model is symmetric with no preferred direction of motion.
This reflects a collective motion of the agents in a random direction given by a Brownian motion.

\begin{figure}[H]
    \centering
    \includegraphics[width=0.87\linewidth]{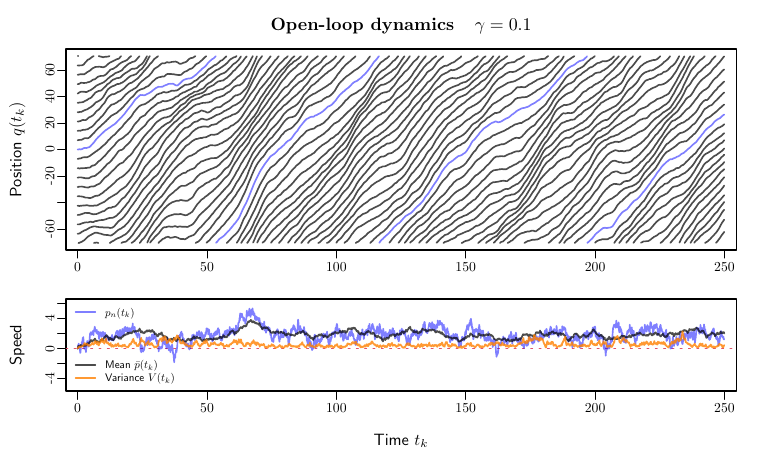}
    \caption{Vehicle trajectories (upper panel) and vehicle mean speed, speed variance and speed of a single vehicle (lower panel) for the uncontrolled system with $\alpha=0.5$, $\beta=1$, $\gamma=0.1$, $u_n(t)=(L/N-\ell)/T$ where $\ell=5$ and $T=1$, $L/N=141/20=7.05$; $\sigma=1$. 
    The system is weakly damped and, although converging, shows evanescent waves propagating at the speed of the vehicles.}
    \label{fig2}
\end{figure}

\begin{figure}[H]
    \centering
    \includegraphics[width=0.87\linewidth]{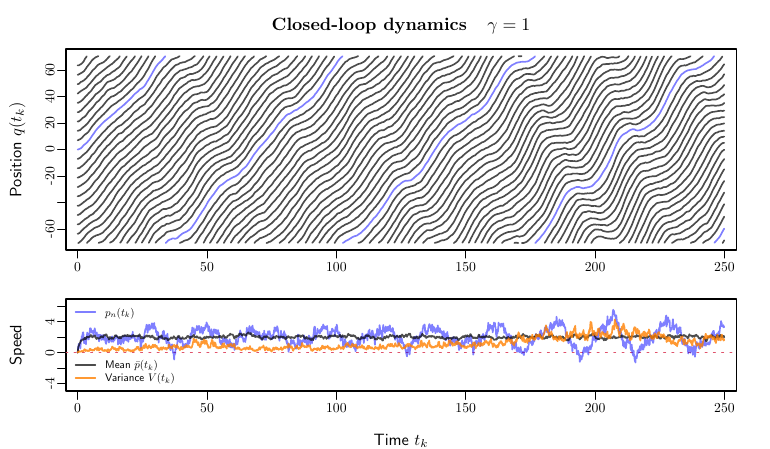}
    \caption{Vehicle trajectories (upper panel) and vehicle mean speed, speed variance and speed of a single vehicle (lower panel) for the uncontrolled system with $\alpha=0.5$, $\beta=1$, $\gamma=1$, $u_n(t)=(\Delta x_n(t)-\ell)/T$ where $\ell=5$ and $T=1$, $L/N=141/20=7.05$; $\sigma=1$. 
    The system is unstable and starts describing stop-and-go waves propagating backward at the characteristic speed $-\ell/T=-5$, forcing the vehicles into stop-and-go dynamics.}
    \label{fig3}
\end{figure}

The simulations of the open-loop and closed-loop systems are shown in Figures~\ref{fig2} and \ref{fig3}, respectively. 
Both models describe wave phenomena. 
The waves propagate at the speed of the control input $x=2.05$ for the open-loop system, especially when, as here with $\gamma=0.1$, the system is weakly damped. 
They appear to be unstable and evanescent. 
They dissipate when the Hamiltonian component $\alpha$, the control rate $\gamma$, or, to a lesser extent, the speed-alignment rate $\beta$ is increased.
In contrast, for the closed-loop system, a single and stable wave moves backward, forcing the vehicles into stop-and-go dynamics. 
Note that $\alpha=0.5$ while $\beta=\gamma=T=1$, so the stability condition $\gamma T+2(\alpha T)^2=1.5 > 2$ does not hold, causing the system to asymptotically explode.

\bibliographystyle{plain}

\begin{thebibliography}{10}

\bibitem{rashad2020twenty}
R.~Rashad, F.~Califano, A.J. van~der Schaft, S.~Stramigioli, 
\textit{Twenty years of distributed port-Hamiltonian systems: a literature review},
IMA J. Math. Contr. Inform. \textbf{37}, 1400 (2020).

\bibitem{knorn2014passivity}
S.~Knorn, A.~Donaire, J.C. Ag{\"u}ero, R.H. Middleton, 
\textit{Passivity-based control for multi-vehicle systems subject to string constraints},
Automatica \textbf{50}, 3224 (2014).

\bibitem{matei2019}
I.~Matei, C.~Mavridis, J.S. Baras, M.~Zhenirovskyy, 
\emph{Inferring particle interaction physical models and their dynamical properties}, in: \emph{2019 IEEE 58th Conference on Decision and Control (CDC)} (IEEE, 2019), pp. 4615--4621.

\bibitem{ehrhardt2024collective}
M.~Ehrhardt, T.~Kruse, A.~Tordeux, 
\textit{The collective dynamics of a stochastic port-Hamiltonian self-driven agent model in one dimension},
ESAIM: Math. Model. Numer. Anal. \textbf{58}, 515 (2024).

\bibitem{ackermann2024stabilisation}
J.~Ackermann, M.~Ehrhardt, T.~Kruse, A.~Tordeux, 
{Stabilisation of stochastic single-file dynamics using port-Hamiltonian systems},
IFAC-PapersOnLine \textbf{58}, 145 (2024).

\bibitem{rudiger2024stability}
B.~R{\"u}diger, A.~Tordeux, B.E. Ugurcan,
\textit{Stability analysis of a stochastic port-Hamiltonian car-following model},
J. Phys. A: Math. Theor. \textbf{57}, 295203 (2024).

\bibitem{Bando1995}
M.~Bando, K.~Hasebe, A.~Nakayama et~al.,
\textit{Dynamical model of traffic congestion and numerical simulation},
Phys. Rev. E \textbf{51}, 1035 (1995).

\bibitem{jiang2001full}
R.~Jiang, Q.~Wu, Z.~Zhu, 
\textit{Full velocity difference model for a car-following theory},
Phys. Rev. E \textbf{64}, 017101 (2001).

\bibitem{Orosz2009}
G.~Orosz, R.E. Wilson, R.~Szalai, G.~St\'ep\'an, 
\textit{Exciting traffic jams: Nonlinear phenomena behind traffic jam formation on highways},
Phys. Rev. E \textbf{80}, 046205 (2009).

\bibitem{Sugiyama2008}
Y.~Sugiyama, M.~Fukui, M.~Kikushi, K.~Hasebe, A.~Nakayama, K.~Nishinari, S.~Tadaki, 
\textit{Traffic jams without bottlenecks. Experimental evidence for the physical mechanism of the formation of a jam},
New J. Phys. \textbf{10}, 033001 (2008).

\bibitem{ISO15822}
{ISO 15622:2018}, 
\emph{{Intelligent transport systems -- Adaptive cruise control systems}} (2018), {International Norm Organisation}.

\bibitem{kesting2008adaptive}
A.~Kesting, M.~Treiber, M.~Sch{\"o}nhof, D.~Helbing,
\textit{Adaptive cruise control design for active congestion avoidance},
Transport. Res. Part C: Emerg. Techn. \textbf{16}, 668 (2008).

\bibitem{Frank1946}
E.~Frank, 
\textit{On the zeros of polynomials with complex coefficients},
Bull. Amer. Math. Soc. \textbf{52}, 144 (1946).

\end{thebibliography}

\end{document}